\documentclass[aap,MSNbibl,seceqn,dvips]{arximspdf}

\doi{10.1214/12-AAP866} 
\volume{23}
\issue{3}
\pubyear{2013}
\firstpage{1129}
\lastpage{1147}

\makeatletter
\newcommand{\eqref}[1]{(\ref{#1})}
\newtheorem{theorem}{Theorem}[section]

\newtheorem{lemma}[theorem]{Lemma}

\newproclaim{example}[theorem]{Example}
\newproclaim{definition}[theorem]{Definition}

\newproclaim{Remark}[theorem]{Remark}

\newcommand{\EE}{\mathbf{E}}
\newcommand{\PP}{\mathbf{P}}
\newcommand{\Var}{\operatorname{Var}}
\newcommand{\Lc}{\mathcal L}
\newcommand{\Lto}{\stackrel{\Lc}{\to}}
\newcommand{\Leq}{\stackrel{\Lc}{=}}
\newcommand{\Lpto}{\stackrel{L^p}{\longrightarrow}}
\newcommand{\Lppto}{\stackrel{L^{p'}}{\longrightarrow}}
\newcommand{\Sc}{\mathcal{S}}
\newcommand{\begt}{
\begin{theorem}}
\newcommand{\ent}{\end{theorem}}
\newcommand{\begl}{
\begin{lemma}}
\newcommand{\enl}{\end{lemma}}
\newcommand{\begr}{
\begin{Remark}}
\newcommand{\enr}{\end{Remark}}
\newcommand{\begx}{
\begin{example}}
\newcommand{\enx}{\end{example}}
\newcommand{\sfrac}[2]{\frac{#1}{#2}}
\newcommand{\eps}{\varepsilon}
\newcommand{\gd}{\delta}
\newcommand{\refS}[1]{Section~\ref{#1}}
\newcommand{\refT}[1]{Theorem~\ref{#1}}
\newcommand{\refL}[1]{Lemma~\ref{#1}}
\newcommand{\refE}[1]{Example~\ref{#1}}
\newcommand{\refR}[1]{Remark~\ref{#1}}
\newcommand{\Holder}{H\"older}
\newcommand{\QuickSort}{\texttt{QuickSort}\ }
\newcommand{\ignore}[1]{}
\makeatother

\begin{document}
\begin{frontmatter}

\title{Distributional convergence for the number of symbol comparisons
used by QuickSort\thanksref{T1}}
\thankstext{T1}{Supported by the Acheson~J.~Duncan Fund for the
Advancement of Research in Statistics.}
\runtitle{Distributional convergence for QuickSort symbols}

\begin{aug}
\author[A]{\fnms{James Allen} \snm{Fill}\corref{}\ead[label=e1]{jimfill@jhu.edu}}
\runauthor{J. A. Fill}
\affiliation{Johns Hopkins University}
\address[A]{Department of Applied Mathematics and Statistics\\
Johns Hopkins University\\
34th and Charles Streets\\
Baltimore, Maryland 21218-2682\\
USA\\
\printead{e1}} 
\end{aug}

\received{\smonth{2} \syear{2012}}
\revised{\smonth{3} \syear{2012}}

%
\begin{abstract}
Most previous studies of the sorting algorithm \texttt{QuickSort} have
used the number of key comparisons as a measure of the cost of
executing the algorithm. Here we suppose that the~$n$ independent and
identically distributed (i.i.d.) keys are each represented as a
sequence of symbols from a probabilistic source and that \texttt{QuickSort}
operates on individual symbols, and we measure the execution
cost as the number of symbol comparisons. Assuming only a mild
``tameness'' condition on the source, we show that there is a limiting
distribution for the number of symbol comparisons after normalization:
first centering by the mean and then dividing by~$n$. Additionally,
under a condition that grows more restrictive as~$p$ increases, we have
convergence of moments of orders~$p$ and smaller. In particular, we
have convergence in distribution and convergence of moments of every
order whenever the source is memoryless, that is, whenever each key is
generated as an infinite string of i.i.d. symbols. This is somewhat
surprising; even for the classical model that each key is an i.i.d.
string of unbiased (``fair'') bits, the mean exhibits periodic
fluctuations of order~$n$.
\end{abstract}

%
\begin{keyword}[class=AMS]
\kwd[Primary ]{60F20}
\kwd[; secondary ]{68W40}.
\end{keyword}
\begin{keyword}
\kwd{QuickSort}
\kwd{limit distribution}
\kwd{$L^p$-convergence}
\kwd{symbol comparisons}
\kwd{natural coupling}
\kwd{probabilistic source}
\kwd{tameness}
\kwd{key comparisons}
\kwd{de-Poissonization}.
\end{keyword}

\end{frontmatter}

\section{Introduction, review of related literature and summary}\label
{Sintroduction}
\subsection{Introduction}\label{SSintroduction}
We consider Hoare's~\cite{h1962} \QuickSort algorithm applied to~$n$
distinct random items (called \textit{keys}) $X_1, \ldots, X_n$,
each represented as a word (i.e., infinite string of symbols such as
bits) from some specified finite or countably infinite alphabet.
We will consider various probabilistic mechanisms
[called \textit{\textup{(}probabilistic\textup{)} sources}] for generating the symbols
within a key, but we will always assume that the keys
themselves are i.i.d. (independent and identically distributed), and we
will later place conditions on the source that rule out the
generation of equal keys.

\QuickSort$(X_1, \ldots, X_n)$ chooses one of the~$n$ keys $X_1, \ldots,
X_n$ (called the ``pivot'') uniformly at random, compares
each of the other keys to it and then proceeds recursively to sort both
the keys smaller than
the pivot and those larger than it.

\textit{Key observation \textup{(}coupling\textup{)}}. Because of the assumption that
the keys are i.i.d., we may take the pivot to be the \textit{first} key
in the sequence, $X_1$. Thus if $X_1, X_2, \ldots$ is an infinite
sequence of keys and $C_n$ is any measure of the cost of sorting~$n$
random keys using any cost function~$c$ (e.g., the number of key
comparisons or the number of symbol comparisons), then we can place all
the random variables~$C_n$ on a common probability space by using $C_n
= c(X_1, \ldots, X_n)$. Notice then that $C_n$ is nondecreasing in~$n$.
We will assume throughout that this natural coupling of the random
variables $C_n$ has been used. The coupling opens up the possibility of
establishing stronger forms of convergence than convergence in
distribution, such as almost sure convergence and convergence in $L^p$,
for suitably normalized $C_n$.

Many authors (Knuth~\cite{k2002v3}, R\'{e}gnier~\cite{r1989},
R\"{o}sler~\cite{r1991}, Knessl and Szpankow\-ski \cite{ks1999},
Fill and Janson~\cite{fj2000smooth,fj2002}, Neininger and
Ruschendorff~\cite{nr2002} and others) have studied $K_n$, the (random)
number of
key comparisons performed by the algorithm. This is an appropriate
measure of the
cost of the algorithm if each key comparison has the same cost. On
the other hand, if keys are represented as words
and comparisons are done by scanning the words from left to right,
comparing the
symbols of matching index one by one, then the cost of comparing two
keys is determined by the number of symbols compared until a difference is
found. We call this number the number of \textit{symbol comparisons} for
the key
comparison, and let $S_n$ denote the total number of symbol comparisons when~$n$
keys are sorted by \texttt{QuickSort}. Symbol-complexity analysis allows
us to compare key-based algorithms
such as {\texttt{QuickSort} with digital algorithms such as those
utilizing digital search trees.

The goal of the present work
is to establish a limiting distribution for the normalized sequence of
random variables $(S_n - \EE S_n) / n$.
Both exact and limiting distributions of $S_n$ will depend on the
source, unlike for~$K_n$.
\subsection{Review of closely related literature (\texttt{QuickSort} and
\texttt{QuickSelect})}\label{SSliterature}
Until now, study of asymptotics for \texttt{QuickSort}'s $S_n$ has been
limited mainly to the expected value $\EE S_n$.
Fill and Janson~\cite{fj2004} were the pioneers in that regard,
obtaining, inter alia, exact and asymptotic expressions
for $\EE S_n$ [consult their Theorem~1.1, and note that the asymptotic
expansion extends through terms of order~$n$ with a
$O(\log n)$ remainder] when the keys are infinite binary strings and
the bits within a key result from i.i.d. fair coin tosses. (We will
refer to this model for key-generation as ``the standard binary
source.'' Equivalently, a key is generated by sampling uniformly from
the unit interval, representing the result in binary notation, and
dropping the leading ``binary point.'') They found that the expected
number of bit comparisons required by {\QuickSort} to sort~$n$ keys is
asymptotically equivalent to $\frac{1}{\ln2} n \ln^2 n$, whereas the
lead-order term of the expected number of \textit{key} comparisons is $2
n \ln n$, smaller by a factor of order $\log n$. Now suppose that $N =
(N(t)\dvtx0 \leq t < \infty)$ is a Poisson process with rate~$1$ and is
independent of the generation of the keys, and let $S(t) := S_{N(t)}$.
The authors also found for each fixed $1 \leq p < \infty$ an upper
bound independent of $t \geq1$ on the $L^p$-norm of
%
%
\begin{equation}
\label{ytbinary}
Y(t) := \frac{S(t) - \EE S(t)}{t},
\end{equation}
(see \cite{fj2004}, Remark 5.1(a), and the corresponding~\cite{fj2010},
Proposition 5.7), leading them to speculate that $Y(t)$ might have a
limiting distribution as $t \to\infty$. We will see that a limiting
distribution does indeed exist, not only for the standard binary source
but for a wide range of sources, as well.

Vall\'{e}e et al.~\cite{vcff2009} greatly extended the scope of~\cite
{fj2004} by establishing for much more general sources both an exact
expression for $\EE S_n$ [consult their Proposition~3 and display~(8)]
and an asymptotic expansion (see their Theorem~1) through terms of
order~$n$ with a $o(n)$ remainder. For the broad class of sources~$\Sc
$ considered, the expected number of symbol comparisons is of lead
order\vspace*{1pt} $\frac{1}{h({\Sc})} n \ln^2 n$, where $h(\Sc)$ is
the entropy of
the source (see their Figure~1 for a definition).

Building on work of Fill and Nakama~\cite{fn2009}, who had in turn
followed closely along the lines of~\cite{fj2004},
Vall\'{e}e et al. \cite{vcff2009} also studied the expected number of
symbol comparisons required by the algorithm
\texttt{QuickSelect}$(n, m)$. This algorithm [aka \texttt{Find}$(n,
m)$], a
close cousin of \texttt{QuickSort} also devised by Hoare~\cite{h1961},
finds a key of specified rank~$m$ from a list of~$n$ keys. The authors
of~\cite{vcff2009} considered the case where $m = \alpha n + o(n)$ for
general $\alpha\in[0, 1]$ [note: we will sometimes refer to \texttt
{QuickQuant}$(n, \alpha)$,
rather than \texttt{QuickSelect}$(n, m)$, in
this case] and a broad class of sources~$\Sc$. They found that the
expected number of symbol comparisons asymptotically has lead term $\rho
_{\Sc}(\alpha) n$, where $\rho_{\Sc}(\alpha)$ is described in their
Figure~1. Unlike in the case of \texttt{QuickSort}, this is only a constant
times larger than the expected number of key comparisons, which is well
known to be asymptotically $\kappa(\alpha) n$ with
\[
\kappa(\alpha) := 2 [1 - \alpha\ln\alpha- (1 - \alpha) \ln(1 - \alpha)].
\]

For either \texttt{QuickSelect} or \texttt{QuickSort}, a deeper probabilistic
analysis of the numbers of key comparisons and symbol comparisons is
obtained by treating entire distributions and not just expectations, in
particular, by finding limiting distributions for suitable
normalizations of these counts and, if possible, establishing
corresponding convergence of moments. Consider
\texttt{QuickQuant}$(n, \alpha)$ first. For both key comparisons and
symbol comparisons a suitable normalization is to divide by~$n$, with
no need to center first. For a literature review on the number of key
comparisons, we refer the reader to
\cite{fn2010}, Section 2.2; the number of symbol comparisons is
discussed next.

Fill and Nakama~\cite{fn2010} (see also~\cite{n2009}) were the first to
establish a limiting distribution for the number of symbol comparisons
for any sorting or searching algorithm. They considered \texttt
{QuickQuant}$(n, \alpha)$ for a broad class of
sources and found a limiting distribution (depending on~$\alpha$, and,
of course, also on the source) for the number $S_n(\alpha)$ of symbol
comparisons (after division by~$n$). It would take us a bit too far
afield to describe the limiting random variable $S(\alpha)$, so we
refer the reader to \cite{fn2010}, Section 3.1, see (3.7), for an
explicit description. In their paper they use the natural coupling
discussed in \refS{SSintroduction} and prove, for each~$\alpha$, that
$S_n(\alpha) / n$ converges to $S(\alpha)$ both (i) almost surely and,
under ever stronger conditions on the source as~$p$ increases, (ii) in
$L^p$. Either conclusion implies convergence in distribution, and (ii)
implies convergence of moments of order $\leq p$. The approach taken
in~\cite{fn2010} is sufficiently general that the authors were able to
unify treatment of key comparisons and symbol comparisons and to
consider various other cost functions (see their Example~2.1).

Now we turn our attention back to \texttt{QuickSort}, the focus of this paper.
Let~$K_n$ (resp., $S_n$) denote the random number of key (resp.,
symbol) comparisons required by \QuickSort to sort a list of~$n$ keys.
We first consider $K_n$, for which we know the following convergence in
law, for some random variable~$T$ (where the immaterial choice of
scaling by
$n + 1$, rather than~$n$, matches with~\cite{r1989}):
%
%
\begin{equation}
\label{kconv}
\frac{K_n -\EE K_n}{n + 1} \Lto T.
\end{equation}
This was proved (i)~by R\'{e}gnier~\cite{r1989}, who used the natural
coupling and martingale techniques to establish convergence both almost
surely and in $L^p$ for every finite~$p$; and (ii)~by
R\"{o}sler~\cite{r1991}, who used the \textit{contraction method} (see
R\"
{o}sler and R\"{u}schendorf~\cite{rr2001} for a general discussion) to
prove convergence in the so-called \textit{minimal $L^p$ metric} for
every finite~$p$ [from which~\eqref{kconv}, with convergence of all
moments, again follows]. An advantage of R\"{o}sler's approach was
identification of the distribution of the limiting~$T$ as the unique
distribution of a zero-mean random variable with finite variance
satisfying the distributional fixed-point equation
%
%
\begin{equation}
\label{fixed}
T \Leq U T + (1-U) T^* + g(U),
\end{equation}
with $g(u) := 1 + 2 u \ln u + 2 (1 - u) \ln(1 - u)$ and where, on the
right-hand side, $T$, $T^*$ and~$U$ are independent random variables,
$T^*$ has the same distribution as~$T$ and $U$ is distributed uniformly
over $(0, 1)$. Later, Fill and Janson~\cite{fj2000fixed} showed that
uniqueness of the zero-mean solution
$\Lc(T)$ to~\eqref{fixed} continues to hold without the assumption of
finite variance, or indeed any other assumption.

\subsection{Summary}\label{SSsummary}
This paper establishes, for a broad class of sources, a limiting
distribution for the number $S_n$ of symbol comparisons for \texttt{QuickSort}.
We tried without success to mimic the approach used in~\cite{fn2010}
for \texttt{QuickQuant}. The approach used in this paper, very broadly
put, is to relate the count $S_n$ of symbol comparisons to various
counts of key comparisons and then rely (heavily) on the result of R\'
{e}gnier~\cite{r1989}. Like Fill and Janson~\cite{fj2004,fj2010}, we
will find it much more convenient to work mainly in continuous time
than in discrete time, but we will also ``de-Poissonize'' our result.
In the continuous-time setting and notation established at~\eqref
{ytbinary} (but without limiting attention to the standard binary
source), we will prove in this paper, assuming that the source is
suitably ``tame'' (in a sense to be made precise), that
%
%
\begin{equation}
\label{main}
Y(t) = \frac{S(t) - \EE S(t)}{t} \Lto Y
\end{equation}
for some random variable~$Y$. Following the lead of~\cite{r1989}
and~\cite{fn2010}, we will use the natural coupling discussed in \refS
{SSintroduction}. Under a mild tameness condition that becomes more
stringent as $p \in[2, \infty)$ increases, we will, in fact, establish
convergence in $L^p$ (see our main \refT{Tmain} for a precise
statement). In particular, for any g-tamed source as defined in Remark
\ref{rem2.3}(a) [e.g., for any (nondegenerate) memoryless source] we have
convergence in $L^p$ for every
finite~$p$. Nondegeneracy of the distribution of~$Y$ is proved by
Bindjeme and Fill~\cite{bfnondeg}; thus the denominator~$t$ used
in~\eqref{main} is not too large to get an interesting limiting distribution.

\textit{Outline of the paper}. After carefully describing in \refS
{SSsource} the probabilistic models used to govern the generation of
keys, reviewing in \refS{SSknown} four known results about the number
of key comparisons we will need in our analysis of symbol comparisons
and listing in \refS{SStools} the other basic probability tools we
will need,
in \refS{Smain} we state and prove our main continuous-time result
about convergence in distribution for the number of symbol comparisons.
We extend the result by de-Poissonization to discrete time in \refS{SdePoi}.

\section{Background and preliminaries}\label{Sprelims}
\subsection{Probabilistic source models for the keys}\label{SSsource}
In this subsection, extracted with only small modifications from~\cite
{fn2010}, we describe what is meant by a probabilistic source (our
model for how the i.i.d. keys are generated) using the terminology and
notation of Vall\'{e}e et al. \cite{vcff2009}.

Let~$\Sigma$ denote a finite totally ordered alphabet (set of symbols),
therefore isomorphic to $\{0, \ldots, r - 1\}$, with\vadjust{\goodbreak} the natural order,
for some finite~$r$; a \textit{word} is then an element of $\Sigma
^{\infty
}$, that is, an infinite sequence (or ``string'') of symbols. We will
follow the customary practice of denoting a word
$w = (w_1, w_2, \ldots)$ more simply by $w_1 w_2 \cdots$.

We will use the word ``prefix'' in two closely related ways. First, the
symbol strings belonging to $\Sigma^k$ are called \textit{prefixes} of
length $k$, and so $\Sigma^* := \bigcup_{0 \leq k < \infty} \Sigma^k$
denotes the set of all prefixes of any nonnegative finite length.
Second, if
$w = w_1 w_2 \cdots$ is a word, then we will call
%
%
\begin{equation}
\label{prefix}
w(k) := w_1 w_2 \cdots w_k \in\Sigma^k
\end{equation}
its \textit{prefix of length~$k$}.

\textit{Lexicographic order} is the linear order (to be denoted in the
strict sense by $\prec$)
on the set of words specified by declaring that $w \prec w'$ if (and
only if) for some $0 \leq k < \infty$ the prefixes of~$w$ and $w'$ of
length~$k$ are equal but $w_{k + 1} < w'_{k + 1}$. Then the
symbol-comparisons cost of determining $w \prec w'$ for such words is
just $k + 1$, the number of symbol comparisons.

A \textit{probabilistic source} is simply a stochastic process $W = W_1
W_2 \cdots$ with state space~$\Sigma$ (endowed with its total $\sigma
$-field) or, equivalently, a random variable~$W$ taking values in
$\Sigma^{\infty}$ (with the product $\sigma$-field). According to
Kolmogorov's consistency criterion (e.g., \cite{c2001}, Theorem
3.3.6), the distributions~$\mu$ of such processes are in one-to-one
correspondence with consistent specifications of finite-dimensional
marginals, that is, of the probabilities
\[
p_w := \mu(\{w_1 \cdots w_k\} \times\Sigma^{\infty}), \qquad w = w_1
w_2 \cdots w_k \in\Sigma^*.
\]
Here the \textit{fundamental probability} $p_w$ is the probability that a
word drawn from~$\mu$ has $w_1 \cdots w_k$ as its length-$k$ prefix.

Because the analysis of \texttt{QuickSort} is significantly more
complicated when its input keys are not all distinct, we will restrict
attention to probabilistic sources with continuous distributions~$\mu$.
Expressed equivalently in terms of fundamental probabilities, our
continuity assumption is that for any $w = w_1 w_2 \cdots\in\Sigma
^{\infty}$ we have
$p_{w(k)} \to0$ as $k \to\infty$, recalling the prefix notation~\eqref
{prefix}.

\begx
\label{Esourcex}
We present a few classical examples of sources. For more
examples, and for further discussion, see \cite{vcff2009}, Section 3.

\begin{longlist}[(a)]
\item[(a)]In computer science jargon, a \textit{memoryless source} is one
with $W_1, W_2, \ldots$ i.i.d. Then the fundamental probabilities $p_w$
have the product form
\[
p_w = p_{w_1} p_{w_2} \cdots p_{w_k},\qquad w = w_1 w_2 \cdots w_k \in
\Sigma^*.
\]

\item[(b)] A \textit{Markov source} is one for which $W_1 W_2 \cdots$
is a
Markov chain.

\item[(c)] An
intermittent source (a model for long-range dependence) over the finite
alphabet $\Sigma= \{0, \ldots, r - 1\}$ is defined by specifying the
conditional distributions $\Lc(W_j | W_1, \ldots, W_{j - 1})$ $(j \geq
2)$ in a way that pays special attention to a particular
symbol~$\underline\sigma$. The source is said to be \textit{intermittent
of exponent $\gamma> 0$ with respect
to~$\underline\sigma$} if $\Lc(W_j | W_1, \ldots, W_{j - 1})$ depends
only on the
maximum value~$k$ such that the last~$k$ symbols in the prefix $W_1
\cdots W_{j - 1}$ are all
$\underline\sigma$ and (i)~is the uniform distribution on~$\Sigma$, if
$k = 0$; and (ii)~if $1 \leq k \leq j -1$, assigns mass $[k / (k +
1)]^{\gamma}$ to~$\underline\sigma$ and distributes the remaining mass
uniformly over the remaining elements of~$\Sigma$.
\end{longlist}
\enx

For our results, the quantity
%
%
\begin{equation}
\label{pikdef}
\pi_k := \max\{p_w\dvtx w \in\Sigma^k\}
\end{equation}
will play an important role, as it did in~\cite{vcff2009}, equation
(7), in connection with the generalized Dirichlet series
$\Pi(s) := \sum_{k \geq0} \pi^{-s}_k$. In particular, it will be
sufficient to obtain $L^p$ convergence in our main result
(\refT{Tmain}) that
%
%
\begin{equation}
\label{picond}
\Pi(- 1 / p) = \sum_{k \geq0} \pi^{1 / p}_k < \infty;
\end{equation}
a sufficient condition for this, in turn, is of course that the source
is $\Pi$-tamed with $\gamma> p$ in the sense of the following definition.

\begin{definition}
Let $0 < \gamma< \infty$
and $0 < A < \infty$. We say that the source is \textit{$\Pi$-tamed \textup{(}with
parameters~$\gamma$ and~$A$\textup{)}} if the sequence $(\pi_k)$ at~\eqref
{pikdef} satisfies
\[
\pi_k \leq A (k + 1)^{- \gamma}\qquad\mbox{for every }k \geq0.
\]
\end{definition}
Observe that a $\Pi$-tamed source is always continuous.

\begin{Remark}\label{rem2.3}
(a)~Many common sources have geometric decrease in $\pi_k$ (call these
``g-tamed'') and so for \textit{any} $\gamma$ are $\Pi$-tamed with
parameters~$\gamma$ and~$A$ for suitably chosen $A \equiv A_{\gamma}$.

For example, a memoryless source satisfies $\pi_k = p^k_{\max}$, where
\[
p_{\max} := \sup_{w \in\Sigma^1} p_w
\]
satisfies $p_{\max} < 1$ except in the highly degenerate case of an
essentially single-symbol alphabet. We also have $\pi_k \leq p^k_{\max
}$ for any Markov source, where now $p_{\max}$ is the supremum of all
one-step transition probabilities, and so such a source is g-tamed
provided $p_{\max} < 1$. Expanding dynamical sources (cf.~\cite{cfv2001})
are also g-tamed.

(b) For an intermittent source as in \refE{Esourcex}, for all
large~$k$ the maximum probability $\pi_k$ is attained by the prefix
$\underline\sigma^k$ and equals
\[
\pi_k = r^{-1} k^{- \gamma}.
\]
Intermittent sources are therefore examples of $\Pi$-tamed sources for
which $\pi_k$ decays at a truly inverse-polynomial rate, not an
exponential rate as in the case of g-tamed sources.
\end{Remark}

\subsection{Known results for the numbers of key comparisons for
\texttt{QuickSort}}\label{SSknown}
In this subsection we review four known \QuickSort key-comparison
results (the first two formulated in discrete time and the next two in
continuous time) that will be useful in proving our main Theorem~\ref
{Tmain}. The first gives exact and asymptotic formulas for the
expected number of key comparisons in discrete time and is extremely
basic and well known. [See, e.g.,~\cite{fj2010}, (2.1)--(2.2).]

\begl
\label{Lfixkeyasy}
Let $K_n$ denote the number of key comparisons required to sort a list
of~$n$ distinct keys. Then
%
%
\begin{eqnarray}\label{fixkeyasy}
\EE K_n &=& 2 (n + 1) H_n - 4 n
\nonumber
\\[-8pt]
\\[-8pt]
\nonumber
&=& 2 n \ln n - (4 - 2 \gamma) n + 2 \ln n + (2 \gamma+ 1) + O(1 / n).
\end{eqnarray}
\enl

The second result [mentioned previously at~\eqref{kconv}] is due to R\'
{e}gnier~\cite{r1989}, who also proved convergence in
$L^p$ for every finite~$p$. Recall the \textit{natural coupling}
discussed in \refS{SSintroduction}.

\begl[\cite{r1989}]
\label{Lregnier}
Under the natural coupling, there exists a random variable~$T$ satisfying
%
%
\begin{equation}
\label{kconvas}
\frac{K_n -\EE K_n}{n + 1} \to T\qquad\mbox{almost surely}.
\end{equation}
\enl

We now shift to continuous time by assuming that the successive keys
are generated at the arrival times of a Poisson process with unit rate.
The number of key comparisons through epoch~$t$ is then $K_{N(t)}$,
which we will abbreviate as $K(t)$; while the sequence $(K_n)$ is
thereby naturally embedded in the continuous-time process, the random
variables $K(n)$ and $K_n$ are not to be confused. We will use such
abbreviations throughout this paper; for example, we will also write
$S_{N(t)}$ as~$S(t)$.

The third result we review is the continuous-time analog of~\refL
{Lfixkeyasy}. Note the difference in constant terms and the much
smaller error term in continuous time.

\begl[(\cite{fj2010}, Lemma 5.1)]
\label{LPoikey}
In the continuous-time setting, the expected number of key comparisons
is given by
\[
\EE K(t) = 2 \int^t_0 (t - y) (e^{-y} - 1 + y) y^{-2} \,dy.
\]
Asymptotically, as $t \to\infty$ we have
%
%
\begin{equation}
\label{Poikeyasy}
\EE K(t) = 2 t \ln t - (4 - 2 \gamma) t + 2 \ln t + (2 \gamma+ 2) +
O(e^{-t} t^{-2}).
\end{equation}
\enl

The fourth result gives bounds on the moments of $K(t)$.
For real $p \in[1, \infty)$, we let $\| W \|_p := ( \EE|W|^p
)^{1 / p}$ denote $L^p$-norm.\vadjust{\goodbreak}

\begl[(\cite{fj2010}, Lemma 5.3)]
\label{LPoikeymoments}
For every real $p \in[1, \infty)$, there exists a constant $c_p <
\infty$ such that
\begin{eqnarray*}
\| K(t) - \EE K(t) \|_p &\leq& c_p t\qquad\mbox{for $t \geq1$},
\\
\| K(t) \|_p &\leq& c_p t^{2 / p}\qquad\mbox{for $t \leq1$}.
\end{eqnarray*}
\enl

In the special case $p = 2$, it follows immediately from \refL
{LPoikeymoments} that
%
%
\begin{equation}
\label{Poikeyvariance}
\Var K(t) \leq c^2_2 t^2 \qquad\mbox{for $0 \leq t < \infty$}.
\end{equation}
%

\subsection{Basic probability tools}\label{SStools}
The following elementary lemma is the basic tool we will use for
$L^p$-convergence. For completeness and the reader's convenience, we
supply a proof.

\begl
\label{Lprob}
Let $Y_k(t)$ be random variables, all defined on a common probability
space, for $k = 0, 1, 2, \ldots$ and $0 \leq t \leq\infty$.
Fix $t_0 \in[0, \infty)$ and $1 \leq p < p' < \infty$ and suppose for
some sequences $(b_k)$ and $(b'_k)$ that:
\begin{longlist}[(ii$'$)]
\item[(i)] for each~$k$ we have $Y_k(t) \to Y_k(\infty)$ almost surely
as $t \to\infty$,
\item[(ii)] for each~$k$ we have $\| Y_k(t) \|_p \leq b_k$ for all $t_0
\leq t < \infty$,
\item[(ii$'$)] for each~$k$ we have $\| Y_k(t) \|_{p'} \leq b'_k <
\infty$ for all $t_0 \leq t < \infty$ and
\item[(iii)] $\sum_{k = 0}^{\infty} b_k < \infty$.
\end{longlist}
Then:
\begin{longlist}[(a)]
\item[(a)] for each $t_0 \leq t \leq\infty$ the series $\sum_{k =
0}^{\infty} Y_k(t)$ converges in $L^p$ to some
random variable~$Y(t)$, and moreover,
\item[(b)] $Y(t) \to Y(\infty)$ in $L^p$ as $t \to\infty$.
\end{longlist}
\enl

\begin{pf}
We assume without loss of generality that $t_0 = 0$. Note that
hypotheses~(ii) and~(ii$'$) extend to $t = \infty$ by Fatou's lemma.
\begin{longlist}[(a)]
\item[(a)] From~(ii) and~(iii) it follows for each $0 \leq t \leq\infty
$ that
the sequence of partial sums $\sum_{k = 0}^K Y_k(t)$, $K = 0, 1, \ldots,$ is a Cauchy sequence in the Banach space $L^p$ and so converges to
some random variable $Y(t)$.

\item[(b)] We first claim for each~$k$ that $Y_k(t) \to Y_k(\infty)$ in $L^p$,
that is, $|Y_k(t) - Y_k(\infty)|^p \to0$ in $L^1$ as $t \to\infty$.
To see this, from~(ii$'$) it follows using \cite{c2001},
Exercise~4.5.8, that $|Y_k(t)|^p$ is uniformly integrable in~$t$, as, therefore,
is $|Y_k(t) - Y_k(\infty)|^p$. Our claim then follows from~(i), since
almost-sure convergence to~$0$ implies convergence in probability
to~$0$, and that together with uniform integrability implies
convergence in $L^1$ (e.g., \cite{c2001}, Theorem 4.5.4).
\end{longlist}

Using the triangle inequality for $L^p$-norm, the claim proved in the
preceding paragraph, and the extended
condition~(ii), we find for any~$K$ that
\[
\limsup_{t \to\infty} \|Y(t) - Y(\infty) \|_p
\leq\limsup_{t \to\infty} \sum_{k = K + 1}^{\infty} \| Y_k(t) -
Y_k(\infty) \|_p
\leq2 \sum_{k = K + 1}^{\infty} b_k.
\]
Now let $K \to\infty$, using~(iii), to complete the proof.\vadjust{\goodbreak}
\end{pf}

Later (\refL{LPoiregnier}) we will transfer \refL{Lregnier} to
continuous time. When we do so, the following result will prove useful.
This \textit{law of the iterated logarithm} (LIL) is well known, and, for
example, can be found for general renewal processes
in~\cite{k1997}, Theorem 12.13.
\begl[(LIL for a Poisson process)]
\label{LLIL}
For a Poisson process~$N$ with unit rate,
%
%
\begin{equation}
\label{LIL}
\PP\biggl( \limsup_{t \to\infty} \frac{N(t) - t}{\sqrt{2 t \ln\ln t}}
= 1, \liminf_{t \to\infty} \frac{N(t) - t}{\sqrt{2 t \ln\ln t}} =
-1 \biggr) = 1.
\end{equation}
\enl
\ignore{
\begin{pf}
By the classical Hartman--Wintner LIL for random walks (e.g., \cite{c2001}, Section~7.5),
%
%
\begin{equation}
\label{LILrw}
\PP( \limsup_{k \to\infty} \frac{N(k) - k}{\sqrt{2 k LL k}} =
1, \liminf_{k \to\infty} \frac{N(k) - k}{\sqrt{2 k LL k}} = -1
) = 1,
\end{equation}
where~$k$ here is an integer parameter. To derive~\eqref{LIL}
from~\eqref{LILrw}, it is clearly sufficient to show that
\[
\PP( \lim_{k \to\infty} \frac{N(k + 1) - N(k)}{\sqrt{2 k LL k}}
= 0 ) = 1.
\]
But this is trivial from the first Borel--Cantelli lemma using the
following tail bound for $N(1) \Leq N(k + 1) - N(k)$ and integer
$b \geq0$:
\[
\hspace{1.6in}\PP(N(1) \geq b) = e^{-1} \sum_{n = b}^{\infty} \frac
{1}{n!} \leq\frac{1}{b!}. \hspace{1.2in}\qed
\]
\noqed
\end{pf}
}

\section{Main results (in continuous time): Convergence in $L^p$ (and
therefore in distribution)}\label{Smain}

The following theorem, which adopts the natural coupling discussed in
\refS{SSintroduction} and utilizes the terminology and notation of
\refS{SSsource} for probabilistic sources, is our main result (for
continuous time).

\begt\label{Tmain}
Consider the continuous-time setting in which
independent and identically distributed
keys are generated from a probabilistic source at the arrival times of
an independent Poisson process~$N$ with unit rate. Let $S(t) =
S_{N(t)}$ denote the number of symbol comparisons required by \QuickSort
to sort the keys generated through epoch~$t$, and let
%
%
\begin{equation}
\label{yt}
Y(t) := \frac{S(t) - \EE S(t)}{t}, \qquad0 < t < \infty.
\end{equation}
Let $p \in[2, \infty)$ and assume that
%
%
\begin{equation}
\label{cond}
\sum_{k = 0}^{\infty} \biggl( \sum_{w \in\Sigma^k} p^2_w \biggr)^{1 /
p} <
\infty.
\end{equation}
Then there exists a random variable $Y$ such that $Y(t) \to Y$ in
$L^p$. Thus $Y(t) \Lto Y$, with convergence of moments of orders
$\leq p$; in particular, \mbox{$\EE Y = 0$}.
\ent

\begr
\label{Rmain}
(a)~Observe that $\sum_{w \in\Sigma^k} p_w = 1$ for each~$k$.
Thus\break
$\sum_{w \in\Sigma^k} p^2_w \leq1$, and condition~\eqref{cond} grows
increasingly stronger as~$p$ increases.\vspace*{-6pt}
\begin{longlist}
\item[(b)] Under the weakest instance $p = 2$ of the assumption~\eqref{cond}
we have $Y(t) \to Y$ in $L^2$, and so $Y(t) \to Y$ in law with
convergence of means and variances. The random variable~$Y$ in \refT
{Tmain} of course does not (more precisely, can be taken not to)
depend on the value of~$p$ considered (because a limit in $L^p$ for
any~$p$ is also a limit in probability, and limits in probability are
almost surely unique).

\item[(c)]~The expected number of symbol comparisons in comparing two
independent keys generated by the given source is
$\sum_{w \in\Sigma^*} p^2_w = \sum_{k = 0}^{\infty} \sum_{w \in\Sigma
^k} p^2_w$. So~\eqref{cond} is certainly sufficient to imply that $\EE
S(t) < \infty$ for every~$t$ [in fact, it follows from calculations to\vadjust{\goodbreak}
be performed in the proof of \refT{Tmain} for $p = 2$ that $\EE
S^2(t) < \infty$ for every~$t$] and that with probability one $S(t) <
\infty$ for all~$t$.

\item[(d)]~The sum on~$w$ in~\eqref{cond} is bounded above by the max-prefix
probability $\pi_k$ defined at~\eqref{pikdef}, and so~\eqref{picond}
(namely, $\sum_k \pi_k^{1 / p} < \infty$) is sufficient for~\eqref
{cond}. Thus from the discussion in \refS{SSsource} we see that \refT
{Tmain} gives $L^p$-convergence for $Y(t)$ for all $\Pi$-tamed sources
with parameter
$\gamma> p$. In particular, for any g-tamed source, such as any
(nondegenerate) memoryless source, we have $Y(t) \to Y$ in $L^p$ for
every $p < \infty$.

\item[(e)]~The standard binary source is a classical example of a periodic
memoryless source; cf.~\cite{vcff2009}---specifically, Definition 3(d),
Theorem 1(ii) and the discussion~(ii) in Section~3. Fill and
Janson~\cite{fj2010}, Proposition 5.4, show explicitly for the standard
binary source that
\[
\EE S(t) =\frac{1}{\ln2} t \ln^2 t - c_1 t \ln t + c_2 t +
\pi_t t + O(\log t)\qquad\mbox{as $t \to\infty$},
\]
where $c_1, c_2$ are explicitly given constants and $\pi_t$ is a
certain periodic function of $\log t$.
Given the periodic term of order~$t$ in the mean for this periodic
source, we find it surprising that \refT{Tmain} nevertheless applies.

\item[(f)]~We wonder (but have not yet considered): Under what
conditions do
we have almost sure convergence in
\refT{Tmain} (or in the discrete-time \refT{Tdiscrete})?
\end{longlist}
\enr

To prepare for the proof of \refT{Tmain}, we ``Poissonize'' \refL{Lregnier}.
\begl
\label{LPoiregnier}
In the continuous-time setting of \refT{Tmain}, let $K(t) =
K_{N(t)}$ denote the number of key comparisons required by
\texttt{QuickSort}. Then for the same random variable~$T$ as in the
discrete-time \refL{Lregnier} we have
\[
\frac{K(t) -\EE K(t)}{t} \to T\qquad\mbox{almost surely as $t \to
\infty$}.
\]
\enl

\begin{pf}
This is routine. According to Lemmas~\ref{Lregnier} and~\ref{Lfixkeyasy},
\[
\frac{K_n -[2 n \ln n - (4 - 2 \gamma) n]}{n + 1} \to T\qquad\mbox
{almost surely as $n \to\infty$}.
\]
Since $N(t) \to\infty$ almost surely as $t \to\infty$, it follows that
\[
\frac{K(t) -[2 N(t) \ln N(t) - (4 - 2 \gamma) N(t)]}{N(t) + 1} \to
T\qquad\mbox{almost surely as $t \to\infty$}.
\]
Using the strong law of large numbers (SLLN) for~$N$ [namely, $N(t) / t
\to1$ almost surely, for which \refL{LLIL} is plenty sufficient], we deduce
\[
\frac{K(t) -[2 N(t) \ln N(t) - (4 - 2 \gamma) t]}{t} \to T\qquad\mbox
{almost surely as $t \to\infty$}.
\]
From the mean value theorem it follows that $|y \ln y - x \ln x| \leq
|y - x| (1 + \ln x + \ln y)$ for $x, y \geq1$. Applying\vadjust{\goodbreak} this with $x =
t$ and $y = N(t)$ and invoking the SLLN and the LIL (\refL{LLIL}), we
find almost surely that for large~$t$ we have
\begin{eqnarray*}
|N(t) \ln N(t) - t \ln t|
&\leq&|N(t) - t| [1 + \ln N(t) + \ln t] \\
&\leq&\sqrt{3 t \ln\ln t} [2
\ln t + 1 + o(1)] \\
&=& O\bigl(\sqrt{t \ln\ln t} \times\ln t\bigr) = o(t),
\end{eqnarray*}
and so
\[
\frac{K(t) -[2 t \ln t - (4 - 2 \gamma) t]}{t} \to T\qquad\mbox{almost
surely as $t \to\infty$}.
\]
The desired result now follows from~\eqref{Poikeyasy} in~\refL{LPoikey}.
\end{pf}

We are now ready for the proof.
\begin{pf*}{Proof of \refT{Tmain}}
We use an idea of Fill and Janson~\cite{fj2010}, Section~5, and
decompose $S(t)$ as $\sum_{k = 0}^{\infty} S_k(t)$, and each $S_k(t)$
further as\break $\sum_{w \in\Sigma^k} S_w(t)$, where for an integer~$k$ and
a prefix $w \in\Sigma^k$ we define (with little possibility of
notational confusion)
\begin{eqnarray*}
S_k(t)
&:= &\mbox{number of comparisons of $(k + 1)$st symbols}, \\
S_w(t)
&:=& \mbox{number of comparisons of $(k + 1)$st symbols
between keys}\\
&& \mbox{with prefix~$w$}.
\end{eqnarray*}
A major advantage of working in continuous time is that,
%
%
\begin{equation}
\label{indep}
\mbox{for each fixed~$k$ and~$t$, the variables $S_w(t)$ with
$w\!\in\!\Sigma^k$ are independent.}\hspace*{-40pt}
\end{equation}
A further key observation, clear after a moment's thought, is this: For
each $w \in\Sigma^*$, as stochastic processes,
%
%
\begin{equation}
\label{replica}\qquad
\bigl(S_w(t)\dvtx t \in[0, \infty)\bigr)\mbox{ is a probabilistic
replica of }
\bigl(K(p_w t)\dvtx t \in[0, \infty)\bigr).
\end{equation}

We define corresponding normalized variables as follows:
\[
Y_k(t) := \frac{S_k(t) - \EE S_k(t)}{t}, \qquad Y_w(t) := \frac{S_w(t)
- \EE S_w(t)}{t},
\]
with the normalized variable $Y(t)$ corresponding to $S(t)$ defined
at~\eqref{yt}. Then
\[
Y(t) = \sum_{k = 0}^{\infty} Y_k(t), \qquad Y_k(t) = \sum_{w \in\Sigma
^k} Y_w(t)\qquad\mbox{($k = 0, 1, \ldots$)}.
\]
To complete the proof of $L^p$-convergence of $Y(t)$ we then need only
to find random variables $Y_k(\infty)$ such that the hypotheses of \refL
{Lprob} are satisfied for some $p' \in(p, \infty)$. [Once we have the
main conclusion of the theorem that $Y(t)$ converges to~$Y$ in $L^p$,
convergence in law with convergence of moments of orders
$\leq p$ follows immediately; in particular, since $\EE Y(t) \equiv0$
and $\EE Y(t) \to\EE Y$, we have $\EE Y = 0$.]

But, for each $w \in\Sigma^*$, the existence of an almost-sure limit,
call it $Y_w(\infty)$, for $Y_w(t)$ as $t \to\infty$ follows
from~\eqref{replica} and \refL{LPoiregnier}; indeed, we see that
$Y_w(\infty)$ has the same distribution as $p_w T$, with~$T$ as in \refL
{LPoiregnier}. Taking the finite sum over $w \in\Sigma^k$, we see
that $Y_k(\infty)$ can be defined as $\sum_{w \in\Sigma^k} Y_w(\infty
)$ to meet hypothesis~(i) of \refL{Lprob}.

To verify the remaining hypotheses we choose $t_0 = 1$ and need to
bound the \mbox{$L^q$-norm} of $Y_k(t)$ for~$k$ a nonnegative integer,
$t \in[1, \infty)$ and $q \in\{p, p'\}$. According to \refL
{LRosenthal} to follow, for any real $q \in[2, \infty)$ there exists
a constant $c'_q$ such that
\[
\|Y_k(t)\|_q \leq c'_q\biggl( \sum_{w \in\Sigma^k} p^2_w \biggr)^{1 / q}
\]
for such~$k$ and~$t$. Thus hypotheses~(ii) and [for any $p' \in(p,
\infty)$] (ii$'$) of \refL{Lprob} hold, and the assumption~\eqref
{cond} implies that~(iii) does as well.
\end{pf*}

\begl
\label{LRosenthal}
Adopt the notation in the above proof of \refT{Tmain}. Then for
every real $q \in[2, \infty)$, there exists a constant $c'_q < \infty$
such that
\[
\|Y_k(t)\|_q \leq c'_q \biggl( \sum_{w \in\Sigma^k} p^2_w \biggr)^{1 / q}
\]
for every nonnegative integer~$k$ and every $t \in[1, \infty)$.
\enl

\begin{pf}
Fix $q \in[2, \infty)$. The first step is to use (as did Fill and
Janson~\cite{fj2010}, proof of Proposition 5.7) Rosenthal's inequality,
relying on the fact [recall~\eqref{indep}] that $S_k(t)$ is the \textit
{independent}
sum of $S_w(t)$ with $w \in\Sigma^k$. According to
Rosenthal's inequality~\cite{rosenthal1970}, Theorem 3 (see also,
e.g., \cite{gut2005}, Theorem~3.9.1), there exists a constant $b_q$ (depending
\textit{only} on~$q$) such that
\begin{eqnarray*}
t^q \| Y_k(t) \|^q_q
&=& \| S_k(t) - \EE S_k(t) \|^q_q \\
&\leq& b_q \max\biggl\{ \sum_{w \in\Sigma^k} \| S_w(t) - \EE S_w(t)
\|^q_q,
\biggl[ \sum_{w \in\Sigma^k} \| S_w(t) - \EE S_w(t) \|^2_2 \biggr]^{q
/ 2} \biggr\}.
\end{eqnarray*}
Utilizing~\eqref{replica} and~\refL{LPoikeymoments} together with the
assumptions $t \geq1$ and $q \geq2$ we therefore find
\begin{eqnarray*}
\| Y_k(t) \|^q_q
&\leq& b_q \max\biggl\{ \sum_{w \in A_k(t)} c^q_q p^q_w
+ \sum_{w \in B_k(t)} (2 c_q)^q p^2_w, \biggl( \sum_{w \in\Sigma^k}
c^2_2 p^2_w \biggr) ^{q / 2} \biggr\} \\
&\leq& b_q \max\biggl\{ (2 c_q)^q \sum_{w \in\Sigma^k} p^2_w, c^q_2
\biggl( \sum_{w \in\Sigma^k} p^2_w \biggr) ^{q / 2} \biggr\} \\
&\leq&(c'_q)^q \sum_{w \in\Sigma^k} p^2_w,
\end{eqnarray*}
where $A_k(t)$ and $B_k(t)$ are the intersections of those $\Sigma^k$
with $\{w\dvtx p_w t \geq1\}$ and $\{w\dvtx p_w t < 1\}$, respectively, and
\[
c'_q := b_q^{1 / q} \max\{ 2 c_q, c_2 \}.
\]
The proof is complete.
\end{pf}

\section{Discrete time} \label{SdePoi}

In this final section we de-Poissonize \refT{Tmain} in order to obtain
an analogous result in discrete time, for which we need to strengthen
the hypothesis slightly.

\begt\label{Tdiscrete}
Let $S_n$ denote the number of symbol comparisons required by
\QuickSort to sort the first~$n$ keys generated.
Let $p \in[2, \infty)$ and assume that for some $p' > p$ we have
%
%
\begin{equation}
\label{condi}
\sum_{k = 0}^{\infty} \biggl( \sum_{w \in\Sigma^k} p^2_w \biggr)^{1 / p'}
< \infty.
\end{equation}
If~$Y$ is the continuous-time limit from \refT{Tmain}, then
%
%
\begin{equation}
\label{maind}
\frac{S_n - \EE S_n}{n} \Lpto Y \qquad\mbox{as $n \to\infty$}.
\end{equation}
In particular, we have convergence in distribution, with convergence of
moments of orders $\leq p$.
\ent

We will derive \refT{Tdiscrete} from \refT{Tmain}, and our proof will
need the following moderate deviation estimate for $N(t)$.

\begl\label{Lmoderate}
For any $0 < \epsilon< 1 / 6$, we have
\[
\PP\bigl(|N(t) - t| \geq t^{(1 / 2) + \epsilon}\bigr) \sim\sqrt{\sfrac
{2}{\pi}}
t^{- \epsilon} \exp\biggl( - \sfrac{1}{2} t^{2 \epsilon} \biggr)\qquad
\mbox{as $t \to\infty$}.
\]
\enl

\begin{pf}
It is well known that the normal approximation gives correct lead-order
asymptotics for right-tail deviations from the mean starting from a
point that is, as here, $o(t^{2/3})$. Thus if~$Z$ is distributed
standard normal, then
\[
\PP\bigl(|N(t) - t| \geq t^{(1 / 2) + \epsilon}\bigr) \sim\PP(|Z| \geq
t^{\epsilon})
\sim\sqrt{\sfrac{2}{\pi}} t^{- \epsilon} \exp\biggl( - \sfrac{1}{2}
t^{2 \epsilon} \biggr)
\]
as claimed.\vadjust{\goodbreak}
\end{pf}

In the following proof, given a sequence of events $(B_n)$, we say that
$B_n$ occurs ``wvlp'' (for ``\textit{with very low probability}'') if
$\PP
(B_n)$ is at most an amount exponentially small in a power of~$n$; we
say that $B_n$ occurs ``wvhp'' (for ``\textit{with very high
probability}'') if the complement $B_n^c$ occurs wvlp.

\begin{pf*}{Proof of \refT{Tdiscrete}}
To prove~\eqref{maind} from the integer-time consequence
\[
\frac{S(n) - \EE S(n)}{n} \Lpto Y
\]
of \refT{Tmain}, it is, of course, sufficient to prove
%
%
\begin{equation}
\label{ss}
\frac{S(n) - S_n}{n} \Lpto0
\end{equation}
and
%
%
\begin{equation}
\label{ee}
\frac{\EE S(n) - \EE S_n}{n} \to0.
\end{equation}
Further, since~\eqref{ee} follows immediately from~\eqref{ss}, it is
sufficient to prove~\eqref{ss}.

To prove~\eqref{ss}, the key is to recall that $S(t) = S_{N(t)}$
where~$N$ is a unit-rate Poisson process independent of $(S_0, S_1,
\ldots)$ and to note that $S_n$ increases with~$n$. Let $0 < \eps< 1 /
3$. Applying \refL{Lmoderate} with $(t, \epsilon)$, there taken to be
$(n + n^{(1/2) + \eps}, \eps/ 2)$, wvhp we have
%
%
\begin{equation}
\label{N+}
N\bigl(n + n^{(1/2) + \eps}\bigr) \geq\bigl(n + n^{(1/2) + \eps}\bigr)
- \bigl(n + n^{(1/2) +
\eps}\bigr)^{{1}/{2} + {1}/{2} \eps} \geq n,
\end{equation}
where the second inequality holds for large enough~$n$. Similarly, wvhp
we have
%
%
\begin{equation}
\label{N-}
N\bigl(n - n^{(1/2) + \eps}\bigr) \leq\bigl(n - n^{(1/2) + \eps}\bigr)
- \bigl(n - n^{(1/2) +
\eps}\bigr)^{{1}/{2} + {1}/{2} \eps} \leq n.
\end{equation}
Because $S_{\cdot} \uparrow$, it follows from \eqref{N+}--\eqref{N-} that
\[
S\bigl(n - n^{(1/2) + \eps}\bigr) \leq S_n \leq S\bigl(n + n^{(1/2) +
\eps}\bigr)\qquad\mbox{wvhp},
\]
and hence, wvhp
\begin{eqnarray*}
|S(n) - S_n|
&\leq&\max\bigl\{S(n) - S\bigl(n - n^{(1/2) + \eps}\bigr), S\bigl(n +
n^{(1/2) + \eps}\bigr) -
S(n)\bigr\} \\
&\leq&\bigl[S(n) - S\bigl(n - n^{(1/2) + \eps}\bigr)\bigr] + \bigl
[S\bigl(n + n^{(1/2) + \eps}\bigr) -
S(n)\bigr] \\
&=& S\bigl(n + n^{(1/2) + \eps}\bigr) - S\bigl(n - n^{(1/2) + \eps}\bigr).
\end{eqnarray*}
So to complete the proof of \refT{Tdiscrete} by proving~\eqref{ss}, it
is sufficient to show that
%
%
\begin{equation}
\label{Sdiff}
\frac{S(n + n^{(1/2) + \eps}) - S(n - n^{(1/2) + \eps})}{n} \Lpto0
\end{equation}
and
%
%
\begin{equation}
\label{wvlp}
\frac{S(n) - S_n}{n} \mathbf{1}(A_n) \Lpto0,\vadjust{\goodbreak}
\end{equation}
where $A_n$ is any event wvlp and $\mathbf{1}(A_n)$ is its indicator. We
prove (a)~\eqref{wvlp} and then (b)~\eqref{Sdiff}.

(a)~To bound the $L^p$-norm of the random variable on the left-hand
side in~\eqref{wvlp}, we use \Holder's inequality
$\| Z_1 Z_2 \|_1 \leq\| Z_1 \|_q \| Z_2 \|_{q'}$ with
\begin{eqnarray*}
Z_1 &=& \biggl| \frac{S(n) - S_n}{n} \biggr|^p,\qquad Z_2 = \mathbf{1}(A_n)^p
= \mathbf{1}(A_n), \\
 q& =& \frac{p'}{p} > 1,\qquad
q' =\frac{p'}{p' - p} > 1;
\end{eqnarray*}
note that $(1 / q) + (1 / q') = 1$, as required. Thus
\begin{eqnarray*}
\biggl\| \frac{S(n) - S_n}{n} \mathbf{1}(A_n)\biggr\|^p_p
&=& \EE\biggl[\biggl| \frac{S(n) - S_n}{n}\biggr| \mathbf{1}(A_n)
\biggr]^p
\\
&\leq&\biggl\| \frac{S(n) - S_n}{n}\biggr\|^p_{p'} \times\PP(A_n)^{1 -
(p / p')}.
\end{eqnarray*}
Because $A_n$ occurs wvlp, it suffices to show that $\|S(n)\|_{p'}$ and
$\|S_n\|_{p'}$ each grow at most polynomially in~$n$.

The first of these two is very easy to handle. Using the
hypothesis~\eqref{condi}, we know from
\refT{Tmain} that
\[
\frac{S(t) - \EE S(t)}{t} \Lppto Y,
\]
and it follows that $\|S(t) - \EE S(t)\|_{p'}$ grows at most linearly
in~$t$ as $t \to\infty$. But from the first sentence of \refR{Rmain}
we see that $\EE S(t)$ grows at most quadratically in~$t$, so by the
triangle inequality $\|S(t)\|_{p'}$ grows at most quadratically in~$t$.

Now we turn our attention to $\|S_n\|_{p'}$. Just as we observed in the
preceding paragraph that $\EE S(t)$ grows at most quadratically
in~$t$, we observe here that
\[
0 \leq S_n \leq\sum_{1 \leq i < j \leq n} C_{i j},
\]
where $C_{i j}$ is the cost of comparing the $i$th and $j$th keys, and
hence (with $C := C_{1 2}$)
\[
\|S_n\|_{p'} \leq\sum_{1 \leq i < j \leq n} \|C_{i j}\|_{p'} = \pmatrix{n
\cr2} \|C\|_{p'}.
\]
So, to conclude that $\|S_n\|_{p'}$ grows at most quadratically in~$n$,
we need only show that $\|C\|_{p'}$ is finite. Indeed, for any $t <
\infty$ we have
\begin{eqnarray*}
\infty
&>& \EE S(t)^{p'} \geq\EE\bigl[ S(t)^{p'} \mathbf{1}\bigl(N(t) \geq
2\bigr)
\bigr] \\
&\geq&\EE\bigl[ C^{p'} \mathbf{1}\bigl(N(t) \geq2\bigr) \bigr] = ( \EE
C^{p'} ) \PP\bigl(N(t) \geq2\bigr)
\end{eqnarray*}
and $\PP(N(t) \geq2) > 0$, so $\EE C^{p'} < \infty$.\vadjust{\goodbreak}

(b)~It remains to establish~\eqref{Sdiff}. From two applications of
\refT{Tmain} it follows quickly that
\begin{eqnarray*}
\frac{S(n + n^{(1/2) + \eps}) - \EE S(n + n^{(1/2) + \eps})}{n} &\Lpto&
Y \quad\mbox{and}\\
\frac{S(n - n^{(1/2) + \eps}) - \EE S(n - n^{(1/2) + \eps})}{n} &\Lpto& Y;
\end{eqnarray*}
thus it suffices to prove
%
%
\begin{equation}
\label{last}
\frac{\EE S(n + n^{(1/2) + \eps}) - \EE S(n - n^{(1/2) + \eps})}{n}
\to0.
\end{equation}

Recall from the proof of \refT{Tmain} that
\[
\EE S(t) = \sum_{k = 0}^{\infty} \sum_{w \in\Sigma^k} \EE K(p_w t)
\]
and from \refL{LPoikey} that we know an explicit formula for $\EE
K(t)$, namely,
\[
\EE K(t) = 2 \int^t_0 (t - y) (e^{-y} - 1 + y) y^{-2} \,dy.
\]
This function and its increasing derivative, call it $d(t)$, are both
easily studied. In particular, $d(t) \sim t$ as
$t \downarrow0$ and $d(t) \sim2 \ln t$ as $t \to\infty$. Hence, for
any $0 < \gd\leq1$ there exists a finite constant
$a_{\gd}$ such that
\[
d(t) \leq a_{\gd} t^{\gd}\qquad\mbox{for all $t \in(0, \infty)$}.
\]
Then, for any $0 < t < u < \infty,$ we have
\begin{eqnarray*}
0
&<& \EE S(u) - \EE S(t) = \sum_{k = 0}^{\infty} \sum_{w \in\Sigma^k}
[ \EE K(p_w u) - \EE K(p_w t) ] \\
&\leq&(u - t) \sum_{k = 0}^{\infty} \sum_{w \in\Sigma^k} p_w d(p_w u)
\leq a_{\gd} b_{\gd} (u - t) u^{\gd}
\end{eqnarray*}
with
\[
b_{\gd} := \sum_{k = 0}^{\infty} \sum_{w \in\Sigma^k} p_w^{1 + \gd}.
\]
Therefore,
\[
\EE S\bigl(n + n^{(1/2) + \eps}\bigr) - \EE S\bigl(n - n^{(1/2) + \eps
}\bigr)
\leq2 a_{\gd} b_{\gd} n^{(1/2) + \eps}\bigl(n + n^{(1/2) + \eps}\bigr
)^{\gd} = o(n)
\]
as desired for~\eqref{last}, provided $\sfrac{1}{2} + \eps+ \gd< 1$
and $b_{\gd} < \infty$. Our proof thus far has been valid for any $0 <
\eps< 1 / 3$, but we now\vspace*{1.5pt} restrict it to $0 < \eps< 1 / 4$ and choose
$\gd= \sfrac{1}{2} - 2 \eps\in
(0, \sfrac{1}{2})$. The proof of \refT{Tdiscrete} will be complete
once we see that~$\eps$ and~$\gd$ can be chosen so that $b_{\gd}$ is finite.

Fix~$k$ and recall that $\sum_{w \in\Sigma^k} p_w = 1$. Let~$V$ be a
random variable with probability mass function $(p_w, w \in\Sigma^k)$,
and let $Z := p^{\gd}_V$. Then
\[
\sum_{w \in\Sigma^k} p_w^{1 + \gd} = \EE Z = \|Z\|_1 \leq\|Z\|_{1 /
\gd} = ( \EE Z^{1 / \gd} )^{\gd}
= (\EE p_V)^{\gd} = \biggl( \sum_{w \in\Sigma^k} p_w^2 \biggr)^{\gd}.
\]
We can arrange for $\gd\geq\sfrac{1}{p'}$ by choosing $0 < \eps\leq
\sfrac{1}{4} - \sfrac{1}{2 p'}$, which is possible because $p' > p \geq
2$. Then
\[
b_{\gd} \leq\sum_{k = 0}^{\infty} \biggl( \sum_{w \in\Sigma^k} p^2_w
\biggr)^{1 / p'} < \infty
\]
by assumption~\eqref{condi}.
\end{pf*}

\section*{Acknowledgment}
We thank Svante Janson for excellent suggestions
that led to improvements to \refT{Tmain}.

%
%

%


\printaddresses


\begin{thebibliography}{24}

\bibitem{bfnondeg}
%
\begin{bmisc}[author]
\bauthor{\bsnm{Bindjeme},~\bfnm{P.}\binits{P.}} \AND
\bauthor{\bsnm{Fill},~\bfnm{J.~A.}\binits{J.~A.}}
(\byear{2012}).
\bhowpublished{The limiting distribution for the number of symbol
comparisons used by
{Q}uick{S}ort is nondegenerate.
Available at \url{http://www.ams.jhu.edu/\textasciitilde fill/}.}
\bptok{imsref}%
\end{bmisc}
%
\endbibitem

\bibitem{c2001}
%
\begin{bbook}[mr]
\bauthor{\bsnm{Chung},~\bfnm{Kai~Lai}\binits{K.~L.}}
(\byear{2001}).
\btitle{A Course in Probability Theory}, \bedition{3rd} ed.
\bpublisher{Academic Press}, \baddress{San Diego, CA}.
\bid{mr={1796326}}
\bptok{imsref}%
\end{bbook}
%
\endbibitem

\bibitem{cfv2001}
%
\begin{barticle}[mr]
\bauthor{\bsnm{Cl{\'e}ment},~\bfnm{J.}\binits{J.}},
\bauthor{\bsnm{Flajolet},~\bfnm{P.}\binits{P.}} \AND
\bauthor{\bsnm{Vall{\'e}e},~\bfnm{B.}\binits{B.}}
(\byear{2001}).
\btitle{Dynamical sources in information theory: A general analysis of trie
structures}.
\bjournal{Algorithmica}
\bvolume{29}
\bpages{307--369}.
\bid{doi={10.1007/BF02679623}, issn={0178-4617}, mr={1887308}}
\bptok{imsref}%
\end{barticle}
%
\endbibitem

\bibitem{fj2000smooth}
%
\begin{bincollection}[mr]
\bauthor{\bsnm{Fill},~\bfnm{James~Allen}\binits{J.~A.}} \AND
\bauthor{\bsnm{Janson},~\bfnm{Svante}\binits{S.}}
(\byear{2000}).
\btitle{Smoothness and decay properties of the limiting {Q}uicksort density
function}.
In \bbooktitle{Mathematics and Computer Science ({V}ersailles, 2000)}
\bpages{53--64}.
\bpublisher{Birkh\"auser}, \baddress{Basel}.
\bid{mr={1798287}}
\bptok{imsref}%
\end{bincollection}
%
\endbibitem

\bibitem{fj2000fixed}
%
\begin{barticle}[mr]
\bauthor{\bsnm{Fill},~\bfnm{James~Allen}\binits{J.~A.}} \AND
\bauthor{\bsnm{Janson},~\bfnm{Svante}\binits{S.}}
(\byear{2000}).
\btitle{A characterization of the set of fixed points of the {Q}uicksort
transformation}.
\bjournal{Electron. Commun. Probab.}
\bvolume{5}
\bpages{77--84 (electronic)}.
\bid{doi={10.1214/ECP.v5-1021}, issn={1083-589X}, mr={1781841}}
\bptok{imsref}%
\end{barticle}
%
\endbibitem

\bibitem{fj2002}
%
\begin{barticle}[mr]
\bauthor{\bsnm{Fill},~\bfnm{James~Allen}\binits{J.~A.}} \AND
\bauthor{\bsnm{Janson},~\bfnm{Svante}\binits{S.}}
(\byear{2002}).
\btitle{Quicksort asymptotics}.
\bjournal{J. Algorithms}
\bvolume{44}
\bpages{4--28}.
\bid{doi={10.1016/S0196-6774(02)00216-X}, issn={0196-6774}, mr={1932675}}
\bptok{imsref}%
\end{barticle}
%
\endbibitem

\bibitem{fj2004}
%
\begin{binproceedings}[mr]
\bauthor{\bsnm{Fill},~\bfnm{James~Allen}\binits{J.~A.}} \AND
\bauthor{\bsnm{Janson},~\bfnm{Svante}\binits{S.}}
(\byear{2004}).
\btitle{The number of bit comparisons used by {Q}uicksort: An average-case
analysis}.
In \bbooktitle{Proceedings of the {F}ifteenth {A}nnual {ACM}-{SIAM} {S}ymposium
on {D}iscrete {A}lgorithms}
\bpages{300--307 (electronic)}.
\bpublisher{ACM}, \baddress{New York}.
\bid{mr={2291065}}
\bptok{imsref}%
\end{binproceedings}
%
\endbibitem

\bibitem{fj2010}
%
\begin{bmisc}[author]
\bauthor{\bsnm{Fill},~\bfnm{J.~A.}\binits{J.~A.}} \AND
\bauthor{\bsnm{Janson},~\bfnm{S.}\binits{S.}}
(\byear{2012}).
\bhowpublished{The number of bit comparisons used by {Q}uicksort: An
average-case
analysis.
Available at \url{http://www.ams.jhu.edu/\textasciitilde fill/}.}
\bptok{imsref}%
\end{bmisc}
%
\endbibitem

\bibitem{fn2009}
%
\begin{barticle}[mr]
\bauthor{\bsnm{Fill},~\bfnm{James~Allen}\binits{J.~A.}} \AND
\bauthor{\bsnm{Nakama},~\bfnm{Tak{\'e}hiko}\binits{T.}}
(\byear{2010}).
\btitle{Analysis of the expected number of bit comparisons required by
{Q}uickselect}.
\bjournal{Algorithmica}
\bvolume{58}
\bpages{730--769}.
\bid{doi={10.1007/s00453-009-9294-3}, issn={0178-4617}, mr={2672478}}
\bptok{imsref}%
\end{barticle}
%
\endbibitem

\bibitem{fn2010}
%
\begin{bmisc}[author]
\bauthor{\bsnm{Fill},~\bfnm{J.~A.}\binits{J.~A.}} \AND
\bauthor{\bsnm{Nakama},~\bfnm{T.}\binits{T.}}
(\byear{2012}).
\bhowpublished{Distributional convergence for the number of symbol
comparisons used by
\texttt{QuickSelect}.
Available at \url{http://www.ams.jhu.edu/\textasciitilde fill/}.}
\bptok{imsref}%
\end{bmisc}
%
\endbibitem

\bibitem{gut2005}
%
\begin{bbook}[mr]
\bauthor{\bsnm{Gut},~\bfnm{Allan}\binits{A.}}
(\byear{2005}).
\btitle{Probability: A Graduate Course}.
\bpublisher{Springer}, \baddress{New York}.
\bid{mr={2125120}}
\bptok{imsref}%
\end{bbook}
%
\endbibitem

\bibitem{h1961}
%
\begin{barticle}[author]
\bauthor{\bsnm{Hoare},~\bfnm{C.~A.~R.}\binits{C.~A.~R.}}
(\byear{1961}).
\btitle{Find (algorithm 65)}.
\bjournal{Communications of the ACM}
\bvolume{4}
\bpages{321--322}.
\bptok{imsref}%
\end{barticle}
%
\endbibitem

\bibitem{h1962}
%
\begin{barticle}[mr]
\bauthor{\bsnm{Hoare},~\bfnm{C.~A.~R.}\binits{C.~A.~R.}}
(\byear{1962}).
\btitle{Quicksort}.
\bjournal{Comput. J.}
\bvolume{5}
\bpages{10--15}.
\bid{issn={0010-4620}, mr={0142216}}
\bptok{imsref}%
\end{barticle}
%
\endbibitem

\bibitem{k1997}
%
\begin{bbook}[mr]
\bauthor{\bsnm{Kallenberg},~\bfnm{Olav}\binits{O.}}
(\byear{1997}).
\btitle{Foundations of Modern Probability}.
\bpublisher{Springer}, \baddress{New York}.
\bid{mr={1464694}}
\bptok{imsref}%
\end{bbook}
%
\endbibitem

\bibitem{ks1999}
%
\begin{barticle}[mr]
\bauthor{\bsnm{Knessl},~\bfnm{Charles}\binits{C.}} \AND
\bauthor{\bsnm{Szpankowski},~\bfnm{Wojciech}\binits{W.}}
(\byear{1999}).
\btitle{Quicksort algorithm again revisited}.
\bjournal{Discrete Math. Theor. Comput. Sci.}
\bvolume{3}
\bpages{43--64 (electronic)}.
\bid{issn={1365-8050}, mr={1695194}}
\bptok{imsref}%
\end{barticle}
%
\endbibitem

\bibitem{k2002v3}
%
\begin{bbook}[author]
\bauthor{\bsnm{Knuth},~\bfnm{D.~E.}\binits{D.~E.}}
(\byear{1998}).
\btitle{The Art of Computer Programming. Volume 3: Sorting and Searching}.
\bpublisher{Addison-Wesley}, \baddress{Reading, MA}.
\bptok{imsref}%
\end{bbook}
%
\endbibitem

\bibitem{n2009}
%
\begin{bmisc}[author]
\bauthor{\bsnm{Nakama},~\bfnm{T.}\binits{T.}}
(\byear{2009}).
\bhowpublished{Analysis of execution costs for {QuickSelect}.
Ph.D. thesis,
Dept. Applied Mathematics and Statistics, Johns Hopkins Univ.
Available at
\url{http://www.ams.jhu.edu/\textasciitilde fill/papers/NakamaDissertation.pdf}.}
\bptok{imsref}%
\end{bmisc}
%
\endbibitem

\bibitem{nr2002}
%
\begin{barticle}[mr]
\bauthor{\bsnm{Neininger},~\bfnm{Ralph}\binits{R.}} \AND
\bauthor{\bsnm{R{\"u}schendorf},~\bfnm{Ludger}\binits{L.}}
(\byear{2002}).
\btitle{Rates of convergence for {Q}uicksort}.
\bjournal{J. Algorithms}
\bvolume{44}
\bpages{51--62}.
\bid{doi={10.1016/S0196-6774(02)00206-7}, issn={0196-6774}, mr={1932677}}
\bptok{imsref}%
\end{barticle}
%
\endbibitem

\bibitem{r1989}
%
\begin{barticle}[mr]
\bauthor{\bsnm{R{\'e}gnier},~\bfnm{Mireille}\binits{M.}}
(\byear{1989}).
\btitle{A limiting distribution for quicksort}.
\bjournal{RAIRO Inform. Th\'eor. Appl.}
\bvolume{23}
\bpages{335--343}.
\bid{issn={0988-3754}, mr={1020478}}
\bptok{imsref}%
\end{barticle}
%
\endbibitem

\bibitem{rosenthal1970}
%
\begin{barticle}[mr]
\bauthor{\bsnm{Rosenthal},~\bfnm{Haskell~P.}\binits{H.~P.}}
(\byear{1970}).
\btitle{On the subspaces of {$L^{p}$} {$(p>2)$} spanned by sequences of
independent random variables}.
\bjournal{Israel J. Math.}
\bvolume{8}
\bpages{273--303}.
\bid{issn={0021-2172}, mr={0271721}}
\bptok{imsref}%
\end{barticle}
%
\endbibitem

\bibitem{r1991}
%
\begin{barticle}[mr]
\bauthor{\bsnm{R{\"o}sler},~\bfnm{Uwe}\binits{U.}}
(\byear{1991}).
\btitle{A limit theorem for ``{Q}uicksort''}.
\bjournal{RAIRO Inform. Th\'eor. Appl.}
\bvolume{25}
\bpages{85--100}.
\bid{issn={0988-3754}, mr={1104413}}
\bptok{imsref}%
\end{barticle}
%
\endbibitem

\bibitem{rr2001}
%
\begin{barticle}[mr]
\bauthor{\bsnm{R{\"o}sler},~\bfnm{U.}\binits{U.}} \AND
\bauthor{\bsnm{R{\"u}schendorf},~\bfnm{L.}\binits{L.}}
(\byear{2001}).
\btitle{The contraction method for recursive algorithms}.
\bjournal{Algorithmica}
\bvolume{29}
\bpages{3--33}.
\bid{doi={10.1007/BF02679611}, issn={0178-4617}, mr={1887296}}
\bptok{imsref}%
\end{barticle}
%
\endbibitem

\bibitem{vcff2009}
%
\begin{bincollection}[mr]
\bauthor{\bsnm{Vall{\'e}e},~\bfnm{Brigitte}\binits{B.}},
\bauthor{\bsnm{Cl{\'e}ment},~\bfnm{Julien}\binits{J.}},
\bauthor{\bsnm{Fill},~\bfnm{James~Allen}\binits{J.~A.}} \AND
\bauthor{\bsnm{Flajolet},~\bfnm{Philippe}\binits{P.}}
(\byear{2009}).
\btitle{The number of symbol comparisons in {Q}uick{S}ort and {Q}uick{S}elect}.
In \bbooktitle{Automata, Languages and Programming. {P}art~{I}}
(\beditor{\binits{S.}\bfnm{S.} \bsnm{Albers}},
\beditor{\binits{A.}\bfnm{A.} \bsnm{Marchetti-Spaccamela}},
\beditor{\binits{Y.}\bfnm{Y.} \bsnm{Matias}},
\beditor{\binits{S.}\bfnm{S.} \bsnm{Nikoletseas}}
\AND
\beditor{\binits{W.}\bfnm{W.} \bsnm{Thomas}}, eds.).
\bseries{Lecture Notes in Computer Science}
\bvolume{5555}
\bpages{750--763}.
\bpublisher{Springer}, \baddress{Berlin}.
\bid{doi={10.1007/978-3-642-02927-1_62}, mr={2544890}}
\bptok{imsref}%
\end{bincollection}
%
\endbibitem

\end{thebibliography}
\end{document}